\documentclass[12pt]{article}
\usepackage[cp850]{inputenc}
\usepackage[T1]{fontenc}\usepackage{amsfonts}
\usepackage{graphicx}
\usepackage{amsmath}
\def\<{\langle}\def\E{\mathbb{E}}\def\R{\mathbb{R}}
\def\>{\rangle}
\def\d{\frac{1}{2}\, }
\def\sign{\mathrm{sign}}
\def\<{\langle}
\def\>{\rangle}
\def\E{\mathbb{E}}\def\D{\mathcal{D}}
\def\R{\mathbb{R}}
\def\inti{\int_{-\infty}^{+\infty}}
\def\intd{\int_{\mathbb{R}^d}}

\title{ The Dirichlet curve of a probability in $\R^d$}
\author{G\'erard Letac\thanks{Universit\'e Paul Sabatier, 118 route de Narbonne, 31062 Toulouse, France}, Mauro Piccioni \thanks{Dipartimento di Matematica, Sapienza Universit\`{a} di Roma, 00185 Rome, Italia}}

\begin{document}\maketitle

\begin{abstract} If $\alpha$ is a probability on $\R^d$ and $t>0,$ consider the Dirichlet random probability $P_t\sim\D(t\alpha) ;$ it is such that for any measurable partition $(A_0,\ldots,A_k)$ of $\R^d$ then $(P_t(A_0),\ldots,P_t(A_k))$ is Dirichlet distributed with parameters $(t\alpha(A_0)\ldots,t\alpha(A_k)).$ If $\int_{\R^d}\log(1+\|x\|)\alpha(dx)<\infty$ the random variable $\int_{\R^d}xP_t(dx)$ of $\R^d$ does exist and we denote by $\mu(t\alpha)$ its distribution. The Dirichlet curve associated to the probability $\alpha$ is the map $t\mapsto \mu(t\alpha).$ It has simple properties like $\lim_{t\searrow 0}\mu(t\alpha)=\alpha$ and $\lim_{t\rightarrow \infty}\mu(t\alpha)=\delta_m$ when $m=\int_{\R^d} x\alpha(dx)$ exists. The present paper shows first that if $m$ exists and if $\psi$ is a  convex function on $\R^d$ then $t\mapsto \int_{\R^d}\psi(x)\mu(t\alpha)(dx)$ is a decreasing function, which means that $t\mapsto \mu(t\alpha)$ is decreasing according to the Strassen convex order of probabilities. The second aim of the paper is to prove a group of results around  the following question: if $\mu(t\alpha)=\mu(s\alpha)$ for some $0\leq s<t$, can we claim that $\mu$ is Cauchy distributed in $\R^d?$

\vspace{4mm}\noindent
\textsc{Keywords:} Dirichlet random probability, Strassen convex order, Cauchy distribution.

\vspace{2mm}
\noindent
\textsc{MSC2010 classification:} 60G57, 62E10.

\end{abstract}

\section{Introduction} If $a_0,\ldots,a_k>0$ and $t=a_0+\cdots+a_k$ recall that the Dirichlet distribution $\D(a_0,\ldots,a_k)$ (as named  by Wilks (1962)) is the law of the random variable $(X_0,\ldots,X_k)$ of $\R^{k+1}$ such that $X_i\geq 0$ for all $i=0,\ldots,k$ and $X_0+\cdots+X_k=1$, with the density of $(X_1,\ldots,X_k)$ equal to
$$\frac{\Gamma(t)}{\Gamma(a_0)\ldots \Gamma(a_k)}(1-x_1-\cdots-x_k)^{a_0-1}x_1^{a_1-1}\ldots x_k^{a_k-1}.$$ 
For $f_0,\ldots,f_k>0$ it satisfies \begin{equation} \label{BD}\E\left(\frac{1}{(f_0X_0+\cdots+f_kX_k)^t}\right)=\frac{1}{f_0^{a_0}\ldots f_k^{a_k}}\end{equation}
See for instance Chamayou and Letac (1991).  By considering moments we can prove  the following weak limits: \begin{eqnarray}\label{LIMD1}\lim_{r\rightarrow \infty}\mathcal{D}(ra_0,\ldots,ra_k)&=&\delta_{(a_0/t,\ldots,a_k/t)}  \\
\label{LIMD2}\lim_{\epsilon \rightarrow 0}\mathcal{D}(\epsilon a_0,\ldots,\epsilon a_k)&=&\sum _{i=0}^k\frac{a_i}{t}\delta_{e_i}\end{eqnarray}
 where $(e_0,\ldots,e_k)$ is the canonical basis of $\R^{k+1}.$

More generally, consider a measured space 
$(\Omega,\frak{A},t\alpha)$ where $\alpha $ is a probability on $\Omega$ and $t>0.$ A quick way to introduce the Dirichlet random probability 
$P_t$ on $\Omega$ associated to  the bounded measure $t\alpha$ follows Sethuraman's stick breaking  method: select independent random variables $B_1,Y_1,\ldots,B_n,Y_n,\ldots$ such that $B_n\sim \alpha$ and $Y_n\sim \beta(1,t)(dy)=t(1-y)^{t-1}\textbf{1}_{(0,1)}(y)dy, $ define $W_1=Y_1$ and for $n>1$ 
$$W_n=Y_n(1-Y_{n-1})\ldots (1-Y_1).$$ It is an easy consequence of the strong law of large numbers that with probability  $1,$ as $N \rightarrow \infty$ then  $\sum_{n=1}^{N}W_n=1-(1-Y_1)\ldots(1- Y_{N-1})\rightarrow 1.$ 
Sethuraman (1994) has proved that the random purely atomic probability $P_t$ on $\Omega$ defined by 
\begin{equation}\label{DD}P_t(dw)=\sum_{n=1}^{\infty}W_n\delta_{B_n}(dw),\end{equation} 
satisfies for any measurable partition $(A_0,\ldots,A_k)$ of $\Omega$
\begin{equation}\label{ED}(P_t(A_0),\ldots,P_t(A_k))\sim \D(t\alpha(A_0)\ldots,t\alpha(A_k)).\end{equation} For this reason the  random probability $P_t$ is said to be a  Dirichlet random probability and its distribution is denoted by $\D(t\alpha).$ One says also that $\alpha$ is the governing probability 
of  $P_t$ and that $t$ is its  intensity.  Of course, $(P_t)_{t\geq 0}$ has a venerable story and the papers by Ferguson (1973), Cifarelli and Regazzini (1990), Diaconis ans Kemperman (1996) and Lijoi and Prunster (2009) are among the  important papers to read on the subject. 

Some simple considerations about $\{\D(t\alpha),\ t>0\}$ are in order. If $f$ is a real bounded measurable function defined on $\Omega$ and if $P_t\sim \D(t\alpha)$ then
 the Fourier transform of the real random variable 
$$X_t(f)=\int_{\Omega} f(w)P_t(dw)=\sum_{n=1}^{\infty} W_n f(B_n)$$  will satisfy for real $s:$ 
\begin{eqnarray}\label{ZD}\lim_{t\rightarrow \infty}\E\left(e^{is\int_{\Omega} f(w)P_t(dw)}\right)&=&e^{is\int_{\Omega} f(w)\alpha(dw)}\\
\label{ID}\lim_{t\searrow 0}\E\left(e^{is\int_{\Omega} f(w)P_t(dw)}\right)&=&\int_{\Omega}e^{is f(w)}\alpha(dw)\end{eqnarray}
If $f$ is taking a finite number of values, this is a reformulation of the statements (\ref{LIMD1}) and (\ref{LIMD2}).
 To show (\ref{ZD}) when $f$ is bounded denote $\alpha(f)=\int_{\Omega}f d\alpha$ for simplicity. Introduce a sequence $g_N$ of functions on $\Omega$ taking a finite number of values such that $\epsilon_N=\sup|g_N-f|\rightarrow_{N\rightarrow \infty}0.$ Then
$$\left|\E(e^{isX_t(f)})-e^{is\alpha(f)}\right|\leq A+B+C$$ where $$A=\left|\E(e^{isX_t(f)})-\E(e^{isX_t(g_N)})\right|,\  B=\left|\E(e^{isX_t(g_N)})-e^{is\alpha(g_N)}\right|,\ C=\left|e^{is\alpha(g_N)}-e^{is\alpha(f)}\right|$$ 
From $|e^{ia}-e^{ib}|\leq |a-b|$ we get $A$ and $C$ are less that $2|s|\epsilon_N.$ Furthermore $\lim_{t\searrow 0}B=0$ since $g_N$ takes  a finite number of values. As a consequence $\limsup_{t\searrow 0}(A+B+C)\leq 2|s|\epsilon_N$ for all $N$ and this proves (\ref{ZD}). The proof of (\ref{ID}) is similar.

Notice that, if we assume  that $\Omega$ is a locally compact separable space, then equality (\ref{ZD}) says that $\lim_{t\rightarrow \infty}\mathcal{D}(t\alpha)=\delta_{\alpha}$
whereas, if we denote by $Q_{\alpha}$ the distribution of the random probability on $\Omega$ defined by $\delta_X$ with $X\sim \alpha$, equality (\ref{ID}) says that $\lim_{t\searrow 0}\mathcal{D}(t\alpha)=Q_{\alpha}$ both in the sense of weak convergence. 

The present paper focuses on the distribution of the random variable $X_t(f)$ when $f$ is neither necessarily non-negative nor bounded, and it can be even valued in $\R^d$ rather than in $\R.$ It is easily seen that  if 
$f:\Omega\rightarrow \R^d$ and $\alpha'$ and $P'_t$ are the respective images by $f$ on $\R^d$ of the probabilities $\alpha$ and $P_t$ on $\Omega$, then $P'_t\sim \D(t\alpha').$ Therefore, in order to study the distribution of $X_t(f)=\int_{\Omega} f(w)P_t(dw)=\int_{\R^d}xP'_t(dx),$ there is no loss of generality in choosing $\Omega=\R^d$ and $f$ equal to the identity.  

The problem of the existence of \begin{equation}\label{MD}X_t=\int_{\R^d}xP_t(dx)=\sum_{n=1}^{\infty}W_nB_n\end{equation} (where now the $B_n$'s are iid,  $\alpha$ distributed in $\R^d$) has been solved by a crucial paper of Feigin and Tweedie (1984) where they prove that 
$\int_{\R^d}\|x\|P_t(dx)<\infty$ almost surely if and only if 
\begin{equation}\label{FG}\int_{\R^d}\log (1+\|x\|)\alpha(dx)<\infty \end{equation}
(actually they did this for $d=1$; the case $d>1$ is easily deduced from it). Let us denote by  $FT_d$  the set of
probabilities $\alpha$ on $\R^d$ such that (\ref{FG}) holds. If $\alpha\in FT_d$ denote by $\mu(t\alpha)$ the distribution in $\R^d$ of $X_t$ defined by (\ref{MD}). We anticipate that  $\mu(t\alpha)\notin FT_d$ in general (see Proposition 6.6 below).

The main character of this paper is  the map $t\mapsto \mu(t\alpha)$ from $(0,\infty)$ to the set of probabilities on $\R^d$. We call this map  the Dirichlet curve associated to the probability $\alpha \in FT_d$ on $\R^d.$ From (\ref{MD}) it is important to observe that if the three random variables 
$X$ (valued in $\R^d),$ $B\sim \alpha$ and $Y\sim \beta(1,t)$ are independent then
\begin{equation}\label{SE}X\sim (1-Y)X+YB\end{equation} if and only if $X\sim X_t.$ This follows from a general result described in Chamayou and Letac (1991) (Proposition 1). It is a useful characterization of $\mu(t\alpha).$

In Proposition 3.4 we  see that $t\mapsto \mu(t\alpha)$ is weakly continuous and that \begin{equation}\label{ZM}\lim_{t\searrow 0}\mu(t\alpha)=\alpha.\end{equation}
Furthermore if \begin{equation}\label{SFG} \int_{\R^d}\|x\|\alpha(dx)<\infty \end{equation}then  $m=\int_{\R^d}x\alpha(dx)$ is well defined and Theorem 3.5 below shows \begin{equation}\label{IM}\lim_{t\rightarrow \infty}\mu(t\alpha)=\delta_m.\end{equation} If $\alpha$ has compact support  these two  facts are  immediate consequences of (\ref{ZD}) and (\ref{ID}). Observe also that (\ref{SFG})  implies through (\ref{MD})  that $\E(X_t)$ exists and is equal to $m,$ for any $t>0.$  Comparing the  behavior of $\mu(t\alpha)$ in the neighbourhood of $0$ and $\infty$, one can make the vague observation that the concentration of $\mu(t\alpha)$ is increasing with $t.$ In order to give  a meaning to  this statement, namely that that for $0\leq s\leq t$ the probability $\mu(t\alpha)$ is more concentrated than $\mu(s\alpha),$ we use the Strassen convex order.  Before stating its definition, let us point out that if $\mu$ is a probability in $\R^d$ having a mean and if $\psi$ is a convex function on $\R^d$ then $\int_{\R^d}\max(0,-\psi(x))\mu(dx) <+\infty.$ This comes from the fact that there exists $a\in \R^d$ and $b\in R$ such that $\psi(x)\geq \<a,x\>+b$  together with the fact that $\mu$ has a mean. As a consequence $\int_{\R^d}\psi(x)\mu(dx)$ makes sense, although it can be possibly $+\infty.$ If $\mu$ and $\nu$ are probabilities on $\R^d$ having means we write $\nu\prec\mu$ if $\int_{\R^d}\psi(x)\nu(dx)\leq \int_{\R^d}\psi(x)\mu(dx)$ for all convex functions $\psi$ on $\R^d.$ Needless to say, this implies that  $\mu$ and $\nu$
have the same mean.

Our main theorem is the following

\vspace{4mm}\noindent \textbf{Theorem 1.1:} If $\int_{\R^d}\|x\|\alpha(dx)<\infty$ then for any  convex function $\psi$ on $\R^d$ and for $0< s\leq t$ we have 
$$\int_{\R^d}\psi(x)\mu(t\alpha)(dx)\leq \int_{\R^d}\psi(x)\mu(s\alpha)(dx)$$
In other terms, $t\mapsto \mu(t\alpha)$ is decreasing for the Strassen convex order on $(0,\infty).$

\vspace{4mm}\noindent We shall comment on this result and we will give examples in Section 2. We will prove it in Section 4, after gathering several properties of $\mu(t\alpha)$ in Section 3.

  Next we suppose that  (\ref{FG}) is fulfilled but not  (\ref{SFG}).
In the asymptotic behavior of $\mu(t\alpha)$ when $t\rightarrow \infty$, Cauchy laws play a crucial role. For $b>0$ and $a\in \R$ denote $w=a+ib$ and consider the Cauchy distribution on $\R$ 
\begin{equation}\label{CD1}c_w(dx)=\frac{1}{\pi}\frac{bdx}{(x-a)^2+b^2}.\end{equation} This notation is borrowed from Letac (1978); it enables us to write  the Fourier transform of $c_w$ in the following  way.  For $s>0$ $$\int_{-\infty}^{\infty}e^{isx}c_w(dx)=e^{isw}.$$ Moreover this formula has a sense for $b=0$, in which case $c_w$ is defined as the Dirac mass $\delta_a.$ 
It is a well know fact due to Yamato (1984) that $\mu(t\alpha)=c_w$ for all $t>0$ when $\alpha=c_w.$ In other terms, the Dirichlet curve governed by  $c_w$ is reduced to a point. If (\ref{SFG}) is not fulfilled, the asymptotic behavior of $\mu(t\alpha)$ is not yet well understood: Theorem 3.5 below shows that the limit points of  $\mu(t\alpha)$ as $t\rightarrow \infty$, are  Cauchy distributions in $\R^d$. In $\R^d,$ what we call a Cauchy distribution is  a probability law such that  all linear forms are  one dimensional Cauchy. They  are carefully studied by Samorodnitsky and Taqqu (1994). We recall in Section 5 some results  about them, particularly the fact that a Cauchy distribution in $\R^d$ has not necessarily a center of symmetry. In Section 6 
we shall study the $\alpha'$s such that  $\mu(t\alpha)=\mu(s\alpha)$ for some $0\leq s<t.$ In many particular cases for $(s,t)$ we will prove that these $\alpha$'s are  Cauchy distributions in $\R^d$.\section{Comments and examples} \textsc{Comments:}
\begin{enumerate} \item The Strassen convex order between probabilities on $\R^d$ has a long story, which is reviewed by  Muller and Stoyan  (2002). Recall that if $\mu$ and $\nu$ are probabilities on $\R^d$ having  a mean, the Strassen theorem (see Strassen (1965)) says that the two following properties are equivalent\begin{itemize}\item for any convex function $\psi$ on $\R^d$ we have 
$\int_{\R^d}\psi(y)\nu(dy)\leq \int_{\R^d}\psi(x)\mu(dx)$ (or $\nu\prec \mu$);\item there exists a probability kernel $K(y,dx)$ from $\R^d$ to itself such that 
$\mu(dx)=\int_{\R^d}\nu(dy)K(y,dx)$ and such that $\int_{\R^d}|x|K(x,dy)$ exists and $\int_{\R^d}xK(y,dx)$ is equal to $y$, $\nu$ almost everywhere;  
(in other terms if $X\sim \mu$ and $Y\sim \nu$ one can find a joint distribution $\nu(dy)K(y,dx)$ for $(X,Y)$ such that $\E(X|Y)=Y).$
\end{itemize}
 \item In practical circumstances, it is difficult to make the kernel $K$ explicit. In particular Theorem 1.1 says $\mu(t\alpha)\prec \mu(s\alpha)$ for $0< s<t$
 but the calculation of $K$ seems to be never possible. 
 
 \item It is useful to know that if $\nu_n\prec \mu_n$ and if $\mu_n$ and $\nu_n$ converge weakly to $\mu$ and $\nu$ respectively,  and if the means of $\mu_n$ and $\nu_n$ converge to the means of $\mu$ and $\nu,$ then $\nu\prec\mu.$ This is Theorem 3.4.6 of Muller and Stoyan (2002). Here is an  application of this fact. With the hypotheses and notations of Theorem 1.1,   we have  $\mu(t\alpha)\prec \alpha$ for any $t>0,$ because of (\ref{ZM}).

 \item If $\mu\prec \nu$ and $\nu\prec \mu$ we have $\mu=\nu.$ To see this in dimension one, use the convex function $\psi_a(x)=(x-a)_+ ,$  getting $\int_{[a,\infty)}(x-a)\mu(dx)=\int_{[a,\infty)}(x-a)\nu(dx)$. Thus $$\int_a^{\infty}\left(\int_{[t,\infty)}\mu(dx)\right)dt=\int_a^{\infty}\left(\int_{[t,\infty)}\nu(dx)\right)dt$$ for any $a$ and $\mu=\nu$ follows. It is easy to pass to higher  dimensions by taking linear forms. 
 \end{enumerate}

 \vspace{4mm}\noindent \textsc{Examples of Strassen convex order:}
 \begin{enumerate}
 \item  A classical example is offered by a  sequence $X_1,\ldots, X_n,\ldots$ of iid random variables of $\R^d$ having a mean. If $\mu_n$ is the distribution of $\overline{X}_n=\frac{1}{n}(X_1+\cdots+X_n)$ then $\mu_n\prec \mu_m$ if $1\leq m\leq n$ since  $\E(\overline{X}_m|\overline{X}_n)=\overline{X}_n.$ For seeing this observe that $j\mapsto \E(X_j|\overline{X}_n)$ does not depend on $j.$ This sequence $(\mu_n)_{n\geq 1}$ presents an analogy with the Dirichlet curve. Indeed, by the weak law of large numbers $\mu_n$ converges weakly to $\delta_{\E(X_1)}$. Moreover, if $X_1\sim c_w$ is Cauchy distributed on $\R$ then $\mu_n\sim c_w,$ for any for any integer $n.$ Furthermore if $\mu_n=\mu_m$ where $m$ is not a rational power of $n$, then $X_1$ is Cauchy or Dirac (see Ramachandran and Rao (1970)).

 \item Suppose that $X\sim \mu,$ $Y\sim \nu$ and $0<U<1$ are independent random variables such that $X\sim (1-U)X+UY$ where $\mu$ and $\nu$ are probabilities on $\R^d$ having  a mean. Then $\mu\prec \nu$, since for any convex function $\psi$, writing $m=\E(U)\in (0,1)$, we obtain 
 $$\E(\psi(X))=\E(\psi((1-U)X+UY))\leq (1-m)\E(\psi(X))+m\E(\psi(Y))\Rightarrow \E(\psi(X))\leq \E(\psi(Y)).$$
 
 \item  To give an explicit example of the above case 2) let us use  the following result   due to Chamayou (2000) (with a different proof). We shall use this proposition in the proof of Theorem 1.1.  
 
 \vspace{4mm}\noindent \textbf{Proposition 2.1:}  Let $0<a<b.$ Let  $X_b\sim \beta(b,b)$, $X_a\sim \beta(a,a)$ and $U\sim \beta(2a,b-a)$ be mutually independent. Then $X_b\sim (1-U)X_b+UX_a.$ 

\vspace{4mm}\noindent \textbf{Proof:} For $|t|<1$ apply (\ref{BD}) to the Dirichlet distribution $(1-U,U)\sim D(b-a,2a)$
and to $f_1=1-tX_b,\ f_2=1-tX_a.$ We get
$$\E(\frac{1}{(1-t(1-U)(X_b+UX_a))^{b+a}})= \E(\frac{1}{(1-tX_b)^{b-a}})\times \E(\frac{1}{(1-tX_a)^{2a}}).$$ 
Now we use the Gauss formula: for $V\sim \beta(B,C-B)$ then 
$$_2F_1(A,B;C;t)=\E\left(\frac{1}{(1-tV)^A}\right).$$ We apply it to $V=X_a$, with $B=a$ and $A=C=2a,$ then to  $V=X_b,$ with $A=b-a$, $B=b$ and $C=2b:$  $$\E(\frac{1}{(1-tX_a)^{2a}})=\frac{1}{(1-t)^a},\ \ \E(\frac{1}{(1-tX_b)^{b\pm a}})=\ _2F_1(b\pm a,b;2b;t).$$ 
Now  we use the Euler formula $$_2F_1(A,B;C,t)=(1-t)^{C-A-B}\ _2F_1(C-A,C-B;C;t).$$
for $A=b-a$, $B=b$ and $C=2b,$ obtaining 
$$\E(\frac{1}{(1-t(1-U)X_b+UX_a))^{b+a}})=\E(\frac{1}{(1-tX_b)^{b+a}})$$ which implies the result. $\square$

As a consequence $\beta(b,b)\prec\beta(a,a)$ if $0<a<b.$ No explicit probability kernel $K(x,dy)$ satisfiying the Strassen characterization for this pair $(\beta(b,b),\beta(a,a))$ is known to us. 
 
 \item Suppose that $\alpha$ is concentrated on $[0,\infty)$ and has a moment of order $n$. Then $G_n(t)=\int_0^{\infty}x^n\mu(t\alpha)(dx)$ exists (see Hjort and Ongaro 2005). Theorem 1.1 implies that $t\mapsto G_n(t)$ is decreasing. Proving directly this fact for small values of $n\geq 2$ is a painful process using classical  inequalities for the moments of $\alpha$, as exemplified by Proposition 3.3 below.

 \end{enumerate}
 
 \vspace{4mm}\noindent \textsc{Examples of Dirichlet curves:}
 \begin{enumerate}
 \item  Bernoulli case: If $\Omega=\R^{d+1}$ and $\alpha=p_0\delta_{e_0}+\cdots+p_d\delta_{e_d}$ where $(e_0,\ldots,e_d)$ is the canonical basis of $\R^{d+1}$ then from (\ref{ED}) we have $P_t=X_0\delta_{e_0}+\cdots+X_d\delta_{e_d}$ where $(X_0,\ldots,X_d)\sim \mathcal{D}(tp_0,\ldots,tp_d).$ This implies that 
$\mu(t\alpha)=\mathcal{D}(tp_0,\ldots,tp_d).$ The fact that in this example we have $\mu(t\alpha)\prec \mu(s\alpha)$ for $0\leq s<t$ is by no means obvious and is a consequence of Theorem 1.1. A particular example is obtained for $d=1$:  the ordinary Bernoulli distribution  $\alpha(dx)=q\delta_0+p\delta_1$ with $p=1-q\in (0,1)$ governs the Dirichlet curve $\mu(t\alpha)=\beta(tp,tq),$ for $t>0.$ Theorem 1.1 shows that, for any $0<a<1$, $$t\mapsto \int_a^1(x-a)\beta(tp,tq)(dx)=\int_0^1(x-a)_+\beta(tp,tq)(dx)$$ is decreasing, a fact that seems quite difficult to prove analytically. For the particular case $p=q=1/2$ Theorem 1.1 is directly obtained by using Proposition 2.1, since $$\mu(t(\d\delta_0+\d\delta_1))=\beta(\frac{t}{2},\frac{t}{2}).$$ 

\item If $\Omega=\R$ and $\alpha(dx)=\beta^{(2)}(\d,\d)(dx)=\frac{1}{\pi}\frac{x^{-\d}}{(1+x)}\textbf{1}_{(0,+\infty)}(x)dx,$
then \begin{equation}\label{CIFM}\mu(t\alpha)(dx)=\beta^{(2)}(t+\d,\d)(dx)=\frac{1}{B(t+\d,\d)}\frac{x^{t-\d}}{(1+x)^{1+t}}\textbf{1}_{(0,\infty)}(x)dx\end{equation}
This is due to Cifarelli and Melilli (1980), later corrected by Hjort and Ongaro (2005). This example has no first moments so Theorem 1.1 cannot be applied to it. However notice that $\lim_{t\rightarrow \infty}\mu(t\alpha)$ does not exist.
\item If $\Omega=\R$ and $\alpha=\beta(\d,\d)$ then $\mu(t\alpha)=\beta(t+\d,t+\d).$ To see this  apply  Lemma 2.1 to the particular case $a=\d$ and $b=t+\d$: the lemma says that if $X\sim \beta(t+\d,t+\d)$, $Y\sim \beta(1,t)$ and $B\sim \beta(\d,\d)$ are independent, then $X\sim (1-Y)X+YB.$ 
From the characterization (\ref{SE}) we get the result. Comparing Example 1 with $d=1$ with the present Example 3 we notice the formula: for $t\geq 1/2$
$$\mu(t\alpha_1)=\beta(\frac{t}{2},\frac{t}{2})=\mu(\frac{t-1}{2}\alpha),$$ with $\alpha_1=(\delta_0+\delta_1)/2$ and $\alpha=\beta(1/2,1/2)$: the curve of $\alpha_1$ contains the curve of $\alpha.$ This is the only example we know in which this happens. 
  \item If $\Omega=\R^2$ and $\alpha$ is the uniform distribution on the circle $\mathbb{U}=\{(x,y)\ ;\ x^2+y^2=1\}$ then $\mu(t\alpha)$ is the distribution of $R_t\Theta$ where $R_t^2\sim \beta(1,t)$ is independent of $\Theta\sim \alpha.$ To see this observe from $(\ref{MD})$ that $\mu(t\alpha)$ must be invariant by rotation since $\alpha$ has this property. Furthermore, the image of  $\alpha$ by the projection $(x,y)\mapsto x$ is also the image of $\beta(\d,\d)$ by $x\mapsto x'=2x-1.$ Using the preceeding example, the image of  $\mu(t\alpha)$ by the projection $(x,y)\mapsto x$ is also the image of $\beta(t+\d,t+\d)$ by $x\mapsto x'=2x-1.$ A slightly tedious calculation leads to the result: for this observe that $X'=R\cos \Theta$ where $\Theta$ is uniform on $(0,2\pi] $ and is independent of $R.$ Therefore if $s>0$ we write $\E(R^{2s})=\E(((X')^2)^s)/\E((\cos^2\Theta)^s).$ Similar examples when $\alpha$ is the uniform distribution on the unit sphere of $\R^d$ with $d>2$ are manageable but they lead to untractable formulas for the distribution of $R_t.$ Explicit calculations about this problem appear in Letac and Piccioni (2014), in the comments following  Theorem 16.

  Already for $d=3$ we are led to deal with the Dirichlet curve of the uniform distribution $\alpha_1$ on $(0,1).$ Diaconis and Kemperman (1994) seem to be the first to have written that
  $$\mu(\alpha_1)= \frac{e}{\pi}\frac{\sin \pi x}{x^{x-1}(1-x)^{-x}}\alpha_1(dx),$$ but $\mu(t\alpha_1)$ for $t\neq 1$ is notoriously complicated, as it can be seen in Lijoi and Prunster (2009). 
\item If $\alpha\in FT_d$, if $X\sim \mu(t\alpha)$ is independent of $U\sim \beta(t,t_0)$ then $XU\sim \mu(t_0\delta_0+t\alpha).$ This remark is due to James (2006). More generally  suppose that $Y=(Y_0,\ldots,Y_n)\sim \D(t_0,\ldots,t_n)$ is independent of $X=(X_0,\ldots,X_n)$, being $X_j\sim \mu(t_j\alpha_j)$ with $\alpha_j \in FT_d,$ for $j=0,\ldots, N.$ Then $$Y_0X_0+\cdots+Y_nX_n\sim \mu (t_0\alpha_0+\cdots+t_n\alpha_n).$$ In particular, for $\alpha_j=\alpha \in FT_d$ for all $j=0,\ldots,N,$ $Y_0X_0+\cdots+Y_nX_n$ still lies on the Dirichlet curve of $\alpha$.
  
\item  This example examines the role of the Cauchy distributions in $\R$ and $\R^d.$  Recall that a Cauchy distribution $c$ in $\R^d$ is a distribution such that if $X\sim c$ then $\<f,X\>$ is Cauchy in $\R$ for any linear form $f$ on $\R^d.$ This means that $\int_{\R^d}e^{is\<f,x\>}c(dx)=e^{isw(f)}$, with $f\mapsto w(f)$ positively homogeneous (that is $w(\lambda f)=\lambda w(f)$ for $\lambda\geq 0)$: the admissible $w$'s will be described in Section 5. If $\alpha$ is a probability on $[0,\infty)$
and if $\rho$ is a probability in $\R^d$ we denote by $\rho\circ \alpha$ the distribution of $XY$ when $X\sim \rho$ and $Y\sim \alpha$ are independent. For $d=1,$ the following invariance principle was obtained by  Yamato (1984) in the particular case $\alpha=\delta_1$ and in general, again for $d=1$   by Hjort and Ongaro (2005):

\vspace{4mm}\noindent \textbf{Proposition 2.2:} If $c$ is Cauchy in $\R^d$ and if $\alpha$
is a probability on $[0,\infty)$ belonging to $ FT_1$ then 
\begin{equation}\label{CHU}\mu(tc\circ \alpha)=c\circ\mu(t\alpha).\end{equation}

 \vspace{4mm}\noindent \textbf{Proof:} The proof is quite easy: since $c$ is Cauchy, then $c\in FT_d$. Furthermore, if $\alpha\in FT_1,$ then $c\circ\alpha\in FT_d$ and $\mu(tc\circ\alpha)$ makes sense. Let $X=(X_n),$ $A=(A_n)$ and $Y=(Y_n)$ be three independent  i.i.d. sequences such that $X_n\sim c$, $A_n\sim \alpha$  and $Y_n\sim \beta(1,t)$ then 
$$\mu(tc\circ\alpha)\sim \sum _{n=1}^{\infty}X_nA_nW_n$$ where $W_1=Y_1$, and $W_n$ denotes $Y_n\prod_{j=1}^{n-1}(1-Y_j)$ as usual. So we have to prove that the latter has the same law as $X_0 \sum _{n=1}^{\infty}A_nW_n$, where $X_0 \sim c$ is independent of everything else.
Recall that the Fourier transform of $c$ is $e^{isw(f)}$, with $w$ positively homogeneous, from which the Fourier transform of $\mu(tc\circ\alpha)$ is obtained as follows:
\begin{eqnarray*}&&\int_{\R^d}e^{is\<f,x\>}\mu(tc\circ\alpha)(dx)=\E\left(\E(e^{is\sum _{n=1}^{\infty}\<f,X_n\>A_nW_n}|A,W)\right)=\E\left(e^{\sum _{n=1}^{\infty}isw_fA_nW_n}\right)\\&&=\int_0^{\infty}e^{isw_f a}\mu(t\alpha)(da)=
\E\left(e^{s\<f,X\>\sum _{n=1}^{\infty}A_nW_n}\right)=
\int_{\R^d}e^{is\<f,x\>}c\circ \mu(t\alpha)(dx).\ \ \square
\end{eqnarray*}
\end{enumerate}

\vspace{4mm}\noindent \textbf{Corollary 2.3:} If $c$ is Cauchy in $\R^d$ then $\mu(tc)=c$ for all $t>0.$ 

\vspace{4mm}\noindent \textbf{Proof:}
Choose $\alpha=\delta_1$ in Proposition 2.2.

\section{Moments and asymptotic properties of the Dirichlet curve} 
 The basic link between $\mu(t\alpha)$ and $\alpha$ is the Proposition 3.1 below, due to Cifarelli and Regazzini (1990).
It is a  considerable extension of (\ref{BD}). For convenience, we give two versions. For a real number $t$ and  a non zero complex number $z$ such that its argument arg $z$ is in $(-\pi,\pi),$ the symbols $\log z$ and $z^t$ mean $\log|z|+i\, \mathrm{arg}(z)$ and $e^{t\log z}.$

\vspace{4mm}\noindent \textbf{Proposition 3.1.} If $\alpha\in FT_1$ then for any real $s$ we have $$\inti \frac{\mu(t\alpha)(dx)}{(1-isx)^t}=e^{-t\inti \log(1-isx)\alpha(dx)}$$
and, for $\Im z>0:$
$$\inti \frac{\mu(t\alpha)(dx)}{(x-z)^t}=e^{-t\inti \log(x-z)\alpha(dx)}$$

\vspace{4mm}\noindent With the methods of Hjort and Ongaro (2005) the next proposition gives informations on the Mellin transform of $\|X\|$ when $X\sim \mu(t\alpha):$

\vspace{4mm}\noindent \textbf{Proposition 3.2.} Let $\alpha\in FT_d.$  Let $X\sim \intd xP(dx),$ where $P\sim \mathcal{D}(t\alpha),$ and let $B\sim \alpha.$ Then for any number $s>0$ we have $$\E(\|X\|^s)<\infty\Leftrightarrow \E(\|B\|^s)<\infty.$$ Under these circumstances, for $d=1$ and if $s=n$ is a positive integer we have  the Hjort-Ongaro formula 
\begin{equation}\label{MOMN}\E(X^n)=\frac {(n-1)!}{(t)_{n-1}}\sum_{k=0}^{n-1}(t)_k\frac{\E(X^k)}{k!}\E(B^{n-k}).\end{equation}
Furthermore if $s\geq 1$ we have $\E(\|X\|^s)\leq \E(\|B\|^s)$ and if $0<s<1$  we have 
\begin{equation}\label{SIS}\frac{\E(\|X\|^s)}{\E(\|B\|^s)}\leq tB(t,s),\ \ \ 
\E(\|B\|^s)\leq \E\left((\inti\|x\|P(dx))^s\right)\end{equation}

\vspace{4mm}\noindent \textbf{Proof.} We prove first the equivalence for $s\geq 1.$ If $X,Y,B$ are independent and $Y\sim \beta_{1,t},$  we have $X\sim (1-Y)X+YB$ from (\ref{SE}). Introduce a random variable $G\sim \gamma_{1+t}$ independent of $X,Y,B$ and observe that
$G'=G(1-Y)\sim \gamma_{t}$ and $G''=GY\sim \gamma_1$ are independent. Therefore \begin{equation}\label{HU}GX\sim G'X+G''B\end{equation} with $X$, $G$, $G'$, $G''$ and $B$ mutually independent. 
For proving part $\Leftarrow,$ we  use (\ref{DD}).
Since $s\geq 1$, one has 
\begin{eqnarray} \nonumber\|X\|^s&\leq &\left(\intd \|x\|P(dx)\right)^s\leq \intd \|x\|^sP(dx)=\sum_{i=1}^{\infty}\|B_i\|^sY_i\prod_{k=1}^{i-1}(1-Y_k),\\ \label{GIULIA} \E(\|X\|^s)&\leq &\E\left(\intd \|x\|^sP(dx)\right)=\E\left(\sum_{i=1}^{\infty}\|B_i\|^sY_i\prod_{k=1}^{i-1}(1-Y_k)\right)=\E(\|B\|^s)<\infty.\end{eqnarray}
For proving part $\Rightarrow,$ let us denote $U=G'X$ and $V=G''B.$ If $\E(\|X\|^{s})<\infty,$ then  $\E(\|U+V\|^s)=\E(G^s)\E(\|X\|^{s})<\infty.$ Denote $C_s(u)=\E(\|u+V\|^s)\leq \infty.$ Since $\E(C_s(U))<\infty$, by Fubini's theorem there exists 
$u_0$ such that $C_s(u_0)<\infty.$   We get  from Minkowski 
$$\E(\|V\|^s)\leq \left(\|u_0\|+(\E(\|V+u_0\|^s))^{1/s}\right)^s<\infty$$  since $V$ is the sum of $V+u_0$ and  the constant $-u_0.$   Since $\E(\|V\|^s)=\E((G'')^s)\E(\|B\|^s)$ we get $\E(\|B\|^s)<\infty$ and part $\Rightarrow$ is proved.  Suppose now that $d=1$  and that $\E(\|B\|^n)<\infty).$ Then (\ref{MOMN}) is easily seen from (\ref{HU}): 
$$(t+1)_n\frac{\E(X^n)}{n!}=\frac{\E(G^nX^n)}{n!}=\sum_{k=0}^n\frac{\E((G')^kX^k)}{k!}\frac{\E((G'')^{n-k}B^{n-k})}{(n-k)!}=\sum_{k=0}^n(t)_k\frac{\E(X^k)}{k!}\E(B^{n-k}).
$$
Subtracting from both sides the $n$-th term of the sum and simplifying one gets the desired expression. Finally assume $0<s<1$ 
and observe that for all $t>0$ we have $(1+t)^s\leq 1+t^s$ (just show that $t\mapsto 1+t^s-(1+t)^s$ is increasing). Together with the triangle inequality, this implies that $\|U+V\|^s\leq \|U\|^s+\|V\|^s$ and therefore by taking expectations
$$(\E(G^s)-\E((G')^s))\E(\|X\|^s)\leq \E((G'')^s)\E(\|B\|^s)$$ which is (\ref{SIS}) since $tB(t,s)=\E((G'')^s)/(\E(G^s)-\E((G')^s)).$ 
For (\ref{SIS}) integrate $x\mapsto \|x\|^s$ with $P(dx)$ defined by (\ref{DD}), use the equality inside (\ref{GIULIA}) and  the following inequality (correct for $0<s<1$)
$$\intd\|x\|^sP(dx)\leq \left(\intd\|x\|P(dx)\right)^s.\ \square$$

\vspace{4mm}\noindent \textbf{Remark.} About the first inequality in (\ref{SIS}) note that $tB(t,s)\geq 1$ for $0<s<1:$ just observe that since $\log \Gamma$ is convex, then $s\mapsto \log tB(t,s)$ is decreasing and zero for $s=1.$

\vspace{4mm}\noindent Next proposition shows that if $\alpha$ is concentrated on $[0,\infty)$ then the first moments of $X_t\sim \mu(t\alpha)$ 
have certain delicate properties (which are probably true for any moment). These properties imply that $t\mapsto \E(X_t^n)$ is decreasing. This fact has been  an incentive for  guessing the statement of Theorem 1.1.

\vspace{4mm}\noindent \textbf{Proposition 3.3.} Let $\alpha$ be a probability on $[0,\infty)$ and  $m_k=\int_0^{\infty}x^{k}\alpha(dx),$ where $k$ is a positive integer.  Let $X_t\sim \mu(t\alpha)$. Suppose that $m_k<\infty$ and consider the function $$c_k(t)=\frac{\E(X_t^k)}{k!}.$$ Then $P_{k-1}(t)=(t+1)_k c_k(t)$ is a polynomial of degree $k-1.$ In particular 
$$P_0(t)=m_1,\ \ P_1(t)=\frac{m_2}{2}+\frac{m_1^2}{2}t,\ \ P_2(t)=\frac{m_3}{3}+\frac{m_1m_2}{2}t+\frac{m_1^3}{6}t^2$$
$$P_3(t)=\frac{m_4}{4}+\left(\frac{m_1m_3}{3}+\frac{m_2^2}{8}\right)t+\frac{m_1^2m_2}{4}t^2+\frac{m_1^4}{24}t^3$$
Finally the polynomial $t\mapsto Q_k(t)= -[(t+1)_k]^2c'_{k+1}(t)$ of degree $2k-1$ has non negative coefficients for $k=1,2,3$. As a consequence the functions $t\mapsto \frac {\E(X_t^n)}{n!}$ are decreasing for $n=2,3,4.$

\vspace{4mm}\noindent \textbf{Proof.} From (\ref{MOMN}) one easily gets $P_0(t)=m_1$ and 
$$P_n(t)=\frac{1}{n+1}m_{n+1}+\frac{t}{n+1}\sum_{k=0}^{n-1}P_k(t)m_{n-k}$$ from which $P_1,P_2,P_3$ are  deduced. One also gets 
$$-[(t+1)_k]^2c'_{k+1}(t)=Q_k(t)=P_k(t)\frac{d}{dt}(1+t)_k-(1+t)_kP'_k(t)$$
%If the $2k-1$ coefficients of $Q_k$ are non negative, clearly $t\mapsto \E(X^{k+1})$ is decreasing. 
The first $Q_k$'s are
$$Q_1(t)=\frac{1}{2}(m_2-m_1^2),\ \ Q_2(t)=(m_3-m_1m_2)+\frac{2}{3}(m_3-m_1^2)t+\frac{m_1}{2}(m_2-m_1^2)t^2,$$
\begin{eqnarray*}Q_3(t)&=&\left(2(m_4-m_1m_3)+\frac{3}{4}(m_4-m_2^2)\right)+3(m_4-m_1^2m_2)t\\&&+\left(\frac{3}{4}(m_4-m_1^2m_2)+\frac{3m_2}{4}(m_2-m_1^2)+2m_1(m_3-m_1m_2)\right)t^2\\&&+\left(\frac{2m_1}{3}(m_3-m_1^3)+\frac{1}{4}(m_2^2-m_1^4)\right)t^3+\frac{m_1^4}{4}(m_2-m_1^2)t^4\end{eqnarray*}
If $B\sim \alpha$ then $m_2-m_1^2=\E((B-m_1)^2)\geq 0,$ $m_4-m_2^2=\E((B^2-m_2)^2)\geq 0$ and 
$$m_3-m_2m_1=\E((B-m_1)^2(B+2m_1))\geq 0, \ m_4-m_3m_1=\E((B-m_1)^2(B^2+m_1B+2m^2_1))\geq 0,$$
$$ m_3- m_1^3=(m_3-m_2m_1)+m_1(m_2-m_1^2)\geq 0,\ \ m_4-m_1^2m_2= (m_4-m_3m_1)+m_1(m_3-m_1m_2)\geq 0.$$
This shows the non negativity of the coefficients of $Q_1,$ $Q_2$ and $Q_3.$ $\square$

\vspace{4mm}\noindent \textbf{Proposition 3.4.} If $\alpha \in FT_d$ then $t\mapsto \mu(t\alpha)$ is weakly continuous on $(0,\infty).$  Furthermore we have  $\lim_{t\searrow 0}\mu(t\alpha)=\alpha.$

\vspace{4mm}\noindent \textbf{Proof.} We fix $t_0>0.$  We consider a sequence $(U_n)_{n\geq 1}$ of iid random variables which are uniform on $(0,1).$ Then $1-U_n^{1/t}\sim \beta(1,t)$. If the $B_n$'s are independent with the same distribution $\alpha$ we consider for $t>0$  and $N$ integer
$$X_{N,t}=\sum_{n=N}^{\infty}(U_1\cdots U_{n-1})^{1/t}(1-U_n^{1/t})B_n.$$
We have $X_t=X_{1,t}\sim \mu(t\alpha).$ Consider $M_{N,t}=\sum_{n=N}^{\infty}(U_1\cdots U_{n-1})^{1/t}\|B_n\|.$ Having $\E(\log (1+ \|B_n\|))$  finite we get $\lim_n \|B_n\|^{1/n}=1$ almost surely.  This comes from  $$\sum_{n=1}^{\infty}\Pr\left(\frac{1}{n}\log(1+\|B_n\|)>\epsilon\right)<\infty$$ and the Borel Cantelli Lemma.  From the law of large numbers we have that $\lim_n\frac{1}{n}\sum_{k=1}^n\log U_k=-1.$ By Cauchy criterion these two remarks imply that $M_{N,t}$ converges. Since $t\mapsto M_{N,t}$ is increasing we conclude that for $0<t\leq t_0$ we have
$$\|X_{N,t}\|\leq M_{N,t}\leq M_{N,t_0}.$$ 
This implies the almost sure uniform convergence of the series $X_t$ on $(0,t_0]$. This implies that $t\mapsto X_t$ is almost surely continuous on $(0,\infty).$ Finally, let us extend the definition of $X_t$ to $t=0$ by $X_0=B_1.$  The above uniform convergence extends to $[0,t_0]$ and $\lim_{t\searrow 0}X_t=B_1$ almost surely. Since almost sure convergence implies weak convergence the proof is complete. $\square$

\vspace{4mm}\noindent \textbf{Theorem  3.5.} If $\intd\|x\|\alpha(dx)<\infty$ and if $m=\intd x\alpha(dx)$ then $\mu(t\alpha)\rightarrow _{t\rightarrow\infty }\delta_m.$  If $\intd\|x\|\alpha(dx)=\infty,$ with $\alpha\in FT_d$,   let $(t_n)$ be a sequence tending to infinity. If  $\mu(t_n\alpha)\rightarrow _{n\rightarrow\infty }\mu$ exists and  is a probability, then $\mu$ is a Cauchy distribution.

\vspace{4mm}\noindent \textbf{Comments.}  In the case  
$\intd\|x\|\alpha(dx)=\infty,$ we have seen in (\ref{CIFM}) that $\lim_{t\rightarrow\infty}\mu(t\alpha)$ may fail to exist. Proposition 2.2 has shown  that if $\alpha$ is the distribution of $M>0$ , if  $C\sim c$ is Cauchy in $\R^d$ and is  independent of $M>0$,  and if $\alpha_1$ is the distribution in $\R^d$ of $MC$, then $\mu(t\alpha_1)$ is the distribution of $X_tC$ where $X_t\sim \mu(t\alpha)$ is independent of $C.$ Now if $\E(M)=m,$ Proposition 2.2 shows that the limit distribution of $X_tC$ is  the Cauchy distribution of $mC.$ This example helped us to guess the second statement of Theorem 3.5. The Dirichlet curve $(\mu(t\alpha))_{t\geq 0}$ is not always tight, as shown by the example (\ref{CIFM}). But even if the Dirichlet curve  is tight, it is not clear that
a limit $\mu(t\alpha)\rightarrow _{t\rightarrow\infty }\mu$ always exists.

\vspace{4mm}\noindent \textbf{Proof of Theorem 3.5.} We assume first that $\intd\|x\|\alpha(dx)<\infty$. It is enough to prove the result for $d=1.$  The idea of the proof is to use Proposition 3.1. For  real $s$ we have 

\begin{equation}\label{LRHS}\inti \frac{\mu(t\alpha)(dx)}{(1-\frac{isx}{t})^t}=e^{-t\inti \log(1-\frac{isx}{t})\alpha(dx)},\end{equation}
We will show  that the left hand side converges to some $\inti e^{isx}\mu(dx)$ and we will show that the right hand side to  converges to $e^{ism}.$
 
For the left hand side of (\ref{LRHS}) we first establish the tightness of the family $\{\mu(t\alpha), t>0\}$.  To see this we consider let $X_t\sim \mu(t\alpha)$ and observe that from Markov inequality and Proposition 3.2 we have for all $t>0:$  
$$\Pr(|X_t|>a)\leq \frac{1}{a}\E(|X_t|)\leq \frac{\E(|B|)}{a}.$$  
Next suppose that for some increasing sequence $(t_n),$ the sequence $\mu(t_{n}\alpha)$ converges weakly to a probability $\mu$ as $n \rightarrow \infty.$
Now we consider \begin{eqnarray*}A(t)&=&\inti \left(\frac{1}{(1-\frac{isx}{t})^{t}}-e^{isx}\right)\mu(t\alpha)(dx)\\B(t)&=&\inti e^{isx}(\mu(t\alpha)(dx)-\mu(dx)).\end{eqnarray*}
The left hand side of (\ref{LRHS}) is $A(t)+B(t)+\inti e^{isx}\mu(dx).$ By Paul L\'evy's theorem the sequence $B(t_{n})$ goes to zero when $k \rightarrow \infty.$

We now show that $\lim_{t\rightarrow \infty}A(t)=0.$ We assume $s\neq 0.$ Let us fix $\epsilon>0$  and $a=\E(|B|)/\epsilon$, and define $$A_0(t)=\int_{|x|\geq a}\left(\frac{1}{(1-\frac{isx}{t})^{t}}-e^{isx}\right)\mu(t\alpha)(dx), \ A_1(t)=A(t)-A_0(t).$$
Since $\int_{|x|\geq a}\mu(t\alpha)
(dx)\leq \epsilon$  and since $|(1-\frac{isx}{t})^{-t}|=(1+\frac{s^2x^2}{t^2})^{-t/2}\leq 1$ we can claim that $A_0(t)\leq 2\epsilon$ for all $t.$ 

Next for $0\leq y<t$ introduce the function 
$$f(t,y)=\frac{1}{(1-\frac{y}{t})^{t}}-e^{y}.$$ This is a non-negative function since $\frac{(t)_n}{t^n}-1\geq 0$ shows 
$f(t,y)=\sum_{n=0}^{\infty}\frac{y^n}{n!}\left(\frac{(t)_n}{t^n}-1\right)>0.$ Furthermore $y\mapsto f(t,y)$ is non-decreasing on $(0,t)$ since 
$\frac{\partial}{\partial y}f(t,y)=\frac{t}{t-y}f(t,y)+\frac{y}{t-y}e^y\leq 0.$
For $-t<sx<t$  we  have 
$$\left|\frac{1}{(1-\frac{isx}{t})^{t}}-e^{isx}\right|=
\left|\sum_{n=0}^{\infty}\frac{(isx)^n}{n!}\left(\frac{(t)_n}{t^n}-1\right)\right|\leq f(t,|sx|)$$
As a consequence, for $t>|sa|$
$$|A_1(t)|\leq \int_{-a}^af(t,|sx|)\mu(t\alpha)(dx)\leq f(t,|sa|)\rightarrow _{t\rightarrow \infty}=0$$
and since the right hand side  goes to $0$ for $t \rightarrow\infty,$ one has $\lim_{t\rightarrow \infty}A(t)=0$. 

For the right hand side of (\ref{LRHS})
we introduce the function $g(t,y)=\frac{t}{2}\log(1+\frac{y^2}{t^2}).$ Now we consider $$-t\inti \log(1-\frac{isx}{t})\alpha(dx)=R(t)+iI(t)$$ where 
$R(t)=-\inti g(t,sx)\alpha(dx)$ and where
\begin{eqnarray*}I(t)&=&-t\inti \mathrm{Arg}(1-\frac{isx}{t})\alpha(dx)=t\inti \arctan\left(\frac{sx}{t}\right)\, \alpha(dx)
\\&=&\inti \left (\int_0^{sx}\frac{t^2dv}{t^2+v^2}\right)\alpha(dx)
\rightarrow_{t\rightarrow \infty}\inti sx\, \alpha(dx)=sm\end{eqnarray*}
(here we have used dominated convergence).  For showing $\lim_{t\rightarrow \infty}R(t)=0$ we fix $\epsilon>0$; we introduce $a>0$ such that $\int_{|sx|>a}|x|\alpha(dx)\leq \epsilon$ and such that $\frac{1}{2}\log(1+y^2)\leq |y|$ if $|y|\geq a.$ Since $y\mapsto g(t,y)$ is increasing we get
$$
|R(t)|=\int_{|sx|\leq a}+\int_{|sx|\geq a} g(t,sx)\alpha(dx)\leq g(t,a)+t\int_{|sx|\geq a} \frac{|sx|}{t}\alpha(dx)\leq g(t,a)+|s|\epsilon
$$ leading to the result since $\lim_{t\rightarrow \infty}g(t,a)=0.$

Finally we have proved that for all probability $\mu$ such that there exists an increasing sequence $(t_n)$ satisfiying $\lim_{n\rightarrow \infty}\mu(t_n\alpha)=\mu$ we have $\inti e^{isx}\mu(dx)=e^{ism}$, that is $\mu=\delta_m.$ This is enough to claim that $\lim_{t\rightarrow \infty}\mu(t\alpha)=\delta_m.$

Let us now assume that $\int_{\R^d}\|x\|\alpha(dx)=\infty$ and  that $\mu(t_n\alpha)\rightarrow _{n\rightarrow \infty }\mu$ exists and  is a probability. We imitate much of the preceeding proof, by starting from (\ref{LRHS}) and proving that $A(t_n)$ and $B(t_n)$ both converge to $0$: the tightness of $(\mu(t_n\alpha))_{t>0}$ is guaranteed by  the existence of $\mu$. Therefore the right hand side of (\ref{LRHS}) has a limit when $n\rightarrow \infty.$ As a consequence, the limit $iw$ of $-\frac{t_n}{s}\inti \log(1-\frac{isx}{t_n})\alpha(dx)$ exists but does not depend on $s>0$ and this implies that the limit of the righthand side of (\ref{LRHS}) is $e^{iws}$, which means that $\mu$ is the one dimensional Cauchy distribution $c_w.$ $\square$

\section{Proof of Theorem 1.1} 

\vspace{4mm}\noindent\textsc{First step.}

The following proposition belongs to folklore (see Hjort and Ongaro (2005) Theorem 1). We  give below  a self-contained proof. In the particular case where $\alpha$ is uniform on the unit sphere of $\R^d$, additional details are given in Section 6 of Letac and Piccioni (2014).

\vspace{4mm}\noindent \textbf{Proposition 4.1:}  If $(W_1,\ldots,W_n)\sim \D(t/n,\ldots,t/n)$  and $B_1,\ldots,B_n$ are independent, with $B_j\sim \alpha\in FT_d$ then the limit distribution of $M_n=W_1B_1+\cdots+W_nB_n$ for $n\rightarrow \infty$ is $\mu(t\alpha).$

\vspace{4mm}\noindent \textbf{Proof:} Let $f\in \R^d$ and $z$ complex with $\Im z>0.$ Then if $W_t\sim \mu(t\alpha)$ we have 
\begin{eqnarray*}
\E\left(\frac{1}{(\<f,M_n\>-z)^t}\right)&=&\E\left(\frac{1}{(\<f,W_1B_1+\cdots+W_nB_n\>-z(W_1+\cdots+W_n))^t}\right)\\&
=&\E\left(\frac{1}{(\<f,B_1\>-z)^{t/n}}\ldots \frac{1}{(\<f,B_n\>-z)^{t/n}}\right)=\left(\E(\frac{1}{(\<f,B_1\>-z)^{t/n}})\right)^n
\end{eqnarray*}
We compute the limit of the last expression as follows. If $z=a+ib$ with $b>0$ write $$e^{U+iV}=\frac{1}{(\<f,B_1\>-a-ib)^{t}}$$ where $U$ and $V$ are real. We have $U\leq -t\log b$ and $V\in (0,\pi).$ Therefore $\E(U)$ makes sense, by allowing $-\infty \leq \E(U).$ Consider now iid random variables  $(U_1,V_1)\ldots, (U_n,V_n)$ with the distribution of $(U,V).$ Then the law of large numbers applies and $\frac{1}{n}(U_1+iV_1+\cdots+U_n+iV_n)$ converges almost surely to $\E(U)+i\E(V)$. 
 Also from $U\leq -t\log b$ we are able to claim that by dominated convergence:
\begin{eqnarray*}\left(\E(\frac{1}{(\<f,B_1\>-z)^{t/n}})\right)^n&=&\E(\exp \frac{1}{n}(U_1+iV_1+\cdots+U_n+iV_n))\\&\rightarrow_{n\rightarrow \infty}&\exp (\E(U)+i\E(V))=e^{-t\E(\log(\<f,B_1\>-z)}\\&=&\E\left(\frac{1}{(\<f,W_t\>-z)^t}\right). \square\end{eqnarray*} 

\vspace{4mm}\noindent\textsc{Second  step.} We want to use Proposition 4.1 in the particular case $n=2^k.$ The reason is that we can realise $\D(t/2^k,\ldots,t/2^k)$ by using products of beta random variables as follows. If $k=1$ and $Z^{t}\sim \beta(\frac{t}{2},\frac{t}{2})$ then
$(W^t_1,W^t_2)=(1-Z^t,Z^t)\sim \D(\frac{t}{2},\frac{t}{2}).$ If $k=2$ and if $Z^t$, $Z^t_0$ and $Z^t_1$ are independent and if $Z^t_i$  are $\beta(\frac{t}{4},\frac{t}{4})$
distributed, then 
\begin{equation}\label{KD}(W^t_1,W^t_2,W^t_3,W^t_4)=((1-Z^t)(1-Z^t_0), \ (1-Z^t)Z^t_0,\ Z^t(1-Z^t_1),\ Z^tZ^t_1)\sim D(\frac{t}{4},\frac{t}{4},\frac{t}{4},\frac{t}{4}).\end{equation}
It is worthwhile to give the details of the proof; taking $f_1,f_2,f_3,f_4>0$ we write 
\begin{eqnarray*}
&&\E\left[\left(f_1W^t_1+f_2W^t_2+f_3W^t_3+f_4W^t_4\right)^{-t}\right]=\\&&\E\left[\left((1-Z^t)(f_1(1-Z^t_0)+f_2Z^t_0)+Z^t(f_3(1-Z^t_1)+f_4Z^t_1)\right)^{-t}\right]=\\&&\E\left[\left((f_1(1-Z^t_0)+f_2Z^t_0)\right)^{-t/2}\right]\times \E\left[\left((f_3(1-Z^t_1)+f_4Z^t_1)\right)^{-t/2}\right]= (f_1f_2f_3f_4)^{-t/4}.
\end{eqnarray*}
More generally the set $\{1,\ldots,2^k\}$ is put in a one to one correspondence $j\mapsto (i_1(j),\ldots,i_k(j))$ with $\{0,1\}^k$ by
$$j=1+\sum_{h=1}^{k}i_h(j)2^{h-1},$$ we introduce for each $h=1,\ldots,k-1$ and each $(i_1,\ldots,i_h)\in \{0,1\}^h$ the random variable $$Z^t_{(i_1,\ldots,i_h)}\sim \beta(\frac{t}{2^{h+1}},\frac{t}{2^{h+1}})$$ in such a way that  these random variables are all independent (and are independent of 
$Z^t$).  We define for $h=1,\ldots,k$ 
\begin{eqnarray*}T^t_{(i_1,\ldots,i_h)}&=&Z^t_{(i_1,\ldots,i_{h-1})}\ \mathrm{if}\ i_h=1,\\
&=&1-Z^t_{(i_1,\ldots,i_{h-1})}\ \mathrm{if}\ i_h=0,\\
W^t_j&=&\prod_{h=1}^{k}T^t_{(i_1(j),\ldots,i_h(j))}. \end{eqnarray*}
One can now prove by induction on $k$ along lines  similar to the case $k=2$ that $(W^t_j)_{j=1}^{2^k}\sim\D(t/2^k,\ldots,t/2^k).$ We skip the details. 

\vspace{4mm}\noindent\textsc{Third   step.} We have seen in the comment following Proposition 2.1 that $0<s<t$ implies that $\beta(t,t)\prec \beta(s,s).$ From Strassen theorem this implies the existence of a probability kernel $K_{s,t}(x,dy)$ on $(0,1)^2$ such that 
$$ K_{s,t}(x,dy)\beta(t,t)(dx)$$ is a joint distribution of $(X,Y)$ with  $X\sim \beta(t,t)$, $Y\sim \beta(s,s)$ and $\E(Y|X)=X.$

\vspace{4mm}\noindent
Next, for fixed $0<s<t$ and each $(i_1,\ldots,i_h)$ with $h=1,\ldots,k-1$  we consider a pair $(Z^s_{(i_1,\ldots,i_h)},Z^t_{(i_1,\ldots,i_h)})$
with respective margins $\beta(\frac{s}{2^{h+1}},\frac{s}{2^{h+1}})$ and $\beta(\frac{t}{2^{h+1}},\frac{t}{2^{h+1}})$ and such that the conditional distribution of the former given the latter is $K_{s/2^{h+1},t/2^{h+1}}.$ Finally all these pairs are mutually independent. Now we create also the $W^s_j$'s from the $Z^s$'s as done in the second step. The important point is now
\begin{eqnarray}\label{KA}&\E(W^s_j|Z^t_{(i_1,\ldots,i_h)}, (i_1,\ldots,i_h)\in \{0,1\}^h, h=0,1,\ldots,k-1)
=\prod_{h=1}^k \E(T^t_{(i_1(j),\ldots,i_h(j))}| Z^t_{(i_1(j),\ldots,i_{h-1}(j))}\\&=\prod_{h=1}^k T^s_{(i_1(j),\ldots,i_h(j))}=W^t_j.\end{eqnarray}
Essentially we are using that if $(X_i,Y_i), i=1,\ldots,n$ are mutually independent pairs of random variables with $X_i$ integrable and $\E(X_i|Y_i)=Y_i$ for $i=1,\ldots,n$, then $\E(\prod_{i=1}^n X_i
|Y_1,\ldots,Y_n)=\prod_{i=1}^n \E(X_i|Y_i)$. From (\ref{KA}) we get \begin{equation}\label{KB}\E(W^s_j|W^t_j)=W^t_j.\end{equation} by using the tower property of conditional expectations: if $\E(X|\mathcal{F})=Y$ then $\E(X|\mathcal{G})=Y$  if $\mathcal{G}\subset \mathcal{F}$ and if $Y$ is $\mathcal{G}$-measurable.

\vspace{4mm}\noindent\textsc{Fourth   step.} For simplicity we continue to omit in the notations $W_j^s$ and $W_j^t$ the fact that these random variables depend on $k.$ Defining like in Proposition 4.1 
$$X^t_k=\sum_{j=1}^{2^k}B_jW_j^t,\ X^s_k=\sum_{j=1}^{2^k}B_jW_j^s$$ we can now claim that from (\ref{KB}) that $$\E(X^s_k|W^t_j,B_j\  \ \forall\  j=1,\ldots,2^k)=X^t_k.$$ Again by the tower property we get  $\E(X^s_k|X^t_k)=X^t_k.$ 
By Strassen theorem this implies that $X^t_k \prec X^s_k$. Furthermore $\E(X^t_k)=\E(X^s_k)=\E(B_1)$ for any integer $k$. By Proposition 4.1 $X^t_k$ and $X^s_k$ converge in law to $\mu(t\alpha)$ and $\mu(s\alpha)$, respectively, as $k \rightarrow \infty$. Moreover these limit distributions keep the same mean vector $\E(B_1)$. The proof of Theorem 1.1 is completed by an application of Comment 3 in Section 2.

\section{Cauchy distributions in $\R^d$} The next problem to deal with is the study of the Dirichlet curve $t\mapsto \mu(t\alpha)$ when $\int_{\R^d}\|x\|\alpha	(dx)=\infty.$ Theorem 3.5 has  essentially shown that if the probability  $\mu(\infty)=\lim_{t\rightarrow \infty}\mu(t\alpha)$ exists then $\mu(\infty)$ is Cauchy in $\R^d$. In Section 6  we will prove various characterizations of the Cauchy distributions related to the Dirichlet curve. These characterisations are linked with   the general conjecture $\mu(t\alpha)=\mu(s\alpha)$ for $t\neq s$ implies that $\alpha$ is Cauchy. To this aim  the present section gives a description of these Cauchy laws.

Recall that we have defined in Section 2  a Cauchy distribution in $\R^d$ as the  distribution of a random vector $X$ such that for each linear form $f$ then $\<f,X\>$ either is Dirac or has a one dimensional Cauchy distribution defined by (\ref{CD1}). In other terms for each $f\in \R^d$ there exists  a complex number $w(f)$ with non negative imaginary part such that $\<f,X\>\sim c_{w(f)}.$
The following proposition   clarifies  the possible $f\mapsto w(f).$ 

\vspace{4mm}\noindent \textbf{Proposition 5.1:} The random variable $X$ in $\R^d$ is Cauchy distributed if and only if there exists $a\in \R^d$ and a positive measure $b(ds)$ on the unit sphere $S$ of $\R^d$ such that $\int_Ss b(ds)=0$  and such that for all $t\in \R^d$ we have $\<f,X\>\sim c_{w(f)}$ with 
 \begin{equation}\label{MC}w(f)=\<a,f\>-\frac{2}{\pi}\int_S \<f,s\>\log |\<f,s\>|b(ds)+i\int_S |\<f,s\>|b(ds)\end{equation}

\vspace{4mm}\noindent \textbf{ Comments:} A remarkable fact about the distribution of  $\<f,X\>$ is that its median  $$\<a,f\>-\frac{2}{\pi}\int_S \<f,s\>\log |\<f,s\>|b(ds)$$ is not a linear form in $f$, which means that the
 distribution of $X$ has not necessarily a center of symmetry. If $b(ds)$ is invariant by $s\mapsto -s$ of course $\int_S \<f,s\>\log |\<f,s\>|b(ds)=0$
 and $a$ is the center of symmetry. 
 
 There are several other definitions of the Cauchy distribution in a Euclidean space in the literature, generally more restrictive that the present one. The most popular is  the distribution of $X$ such that $\E\left(e^{\<t,X\>}\right)=e^{-\|t\|}$ and its affine deformations. For such an  $X$ we have $w(f)=i\|f\|$ and $b(ds)=CU(ds)$ where $U(ds)$ is the uniform probability on the unit sphere $S$ and $C=\sqrt{\pi}\Gamma((d+1)/2)/\Gamma(d/2).$

 For an example of a Cauchy distribution in $\R^2$ without center of symmetry one can consider $b=\delta_1+\delta_{j}+\delta_{j^2}$ where $S$ is identified with the unit circle of the complex plane and where $j$ and $j^2$ are the complex cubic roots of the unity. It satisfies $\int_Ssb(ds)=0.$ If $f=e^{i\theta}$ and if $g(\theta)=-\frac{2}{\pi}\cos \theta\log |\cos \theta|$ then the median of $\<f,X\>$ is 
 $$r(\theta)=g(\theta)+g(\theta-\frac{2\pi}{3})+g(\theta+\frac{2\pi}{3}).$$
 and  $\theta\mapsto r(\theta)e^{i\theta}$ is the equation of a nice trefoil curve.

\vspace{4mm}\noindent \textbf{Proof:} We follow the definitions of Samorodnitsky and Taqqu (1994) chapter 2 and use their results. Since $\<f,X\>$ either is Dirac or has a one dimensional Cauchy distribution, this implies that $X$ is  1 stable. Therefore (Theorem 2.3.1) there exists $a\in R^d$ and a positive measure $b$ on $S$ such that $\E(e^{i\<t,X\>})=e^{\psi(t)}$ where $$\psi(t)=i\<t,a\>-\int_S|\<t,s\>|\left(1+i\frac{2}{\pi}\sign (\<t,s\>)\log |\<t,s\>|\right)b(ds)$$Furthermore, $X$ is strictly 1 stable and  from Theorem 2.4.1 we have $\int_Ssb(ds)=0.$ Writing $t=rf$ with $r>0$ in the above formula, we get  $\E(e^{ir\<f,X\>})=e^{irw(f)},$ where $w(f)$ is given by (\ref{MC}). $\square$

\section{Cauchy distribution and Dirichlet curve} 

All along this section we exploit the properties of the Stieltjes transform of a probability $\alpha$ on $\mathbb{R}$, namely the function, defined for all complex numbers $z$ with $\Im z>0$ by $y(z)=\inti \frac{\alpha(dw)}{w-z}.$ Recall that the Stieltjes transform of the Cauchy distribution  $c_{w}$ with $w=a+ib\in H_+$ and $\overline{w}=a-ib$ is $$\inti \frac{c_{w}(dt)}{t-z}=\frac{1}{\overline{w}-z}$$ 
To start with, for any positive integer $k$ we have $y^{(k)}(z)=k!\inti \frac{\alpha(dw)}{(w-z)^{k+1}}.$

\vspace{4mm}\noindent\textbf{Proposition 6.1.} Let $\alpha\in FT_1$ and let $y$ be its Stieltjes transform. Then $\mu(n\alpha)=\alpha$ if and only if \begin{equation}\label{EQD1}ny(z)y^{(n-1)}(z)=y^{(n)}(z)\end{equation} In particular for $n=1$ and $n=2$ this implies that $\alpha$ is Cauchy or Dirac. If $\alpha\in FT_d$ again $\mu(\alpha)=\alpha$ or $\mu(2\alpha)=\alpha$ if and only if $\alpha$ is Cauchy in $\R^d.$ 

\vspace{4mm}\noindent \textbf{Proof.} Suppose $d=1$ and  use Proposition 3.1.  If $\mu(n\alpha)=\alpha \in FT_1$ we can write with $g(z)=-\inti \log(w-z)\alpha(dw):$ 
$$\inti \frac{\alpha(dw)}{(w-z)^n}=e^{ng(z)}.$$ Both sides are analytic functions on the half plane $H^+=\{z\in \mathbb{C}: \Im z>0\}$.  Deriving in $z$ and using $y=g'$ we get 
$$n\inti \frac{\alpha(dw)}{(w-z)^{n+1}}=ne^{ng(z)}g'(z)=ny(z)\inti \frac{\alpha(dw)}{(w-z)^n},$$ 
from which (\ref{EQD1}) is immediate. Conversely, from (\ref{EQD1}) we write 
$$ny(z)=ng'(z)=\frac{y^{(n)}(z)}{y^{(n-1)}(z)}$$ and we get that
$ y^{(n-1)}$ is proportional to $e^{ng}.$ Since, up to a muliplicative  constant, the left hand side is equal to $\inti \frac{\alpha(dw)}{(w-z)^n}$, we get for some constant $C$
$$\inti \frac{\alpha(dw)}{(w-z)^n}=Ce^{ng(z)}.$$ To see that $C=1$ we use the fact that $\alpha$ has mass $1$ and we replace $z$ by $ri$ with $r>0$ in the equality. We get $$\inti r^n\frac{\alpha(dw)}{(w-ri)^n}=Ce^{n(g(ri)+\log r)}.$$ Now $\lim_{r\rightarrow\infty}\inti r^n\frac{\alpha(dw)}{(w-ri)^n}=i^n.$ Also $$g(ri)+\log r=-\inti\log(\frac{w}{r}-i)\alpha(dw)\rightarrow_{r\rightarrow\infty}\log(-i)=-\frac{\pi}{2}i$$ and therefore $\lim_{r\rightarrow\infty}e^{n(g(ri)+\log r)}=e^{n\frac{\pi}{2}i}=i^n$ which implies $C=1.$

As far as the second statement is concerned, for $n=1$ this is a result due to Lijoi and Regazzini (2004). Our proof is shorter, since the general solution of the differential equation $y'(z)=y^2(z),$ corresponding to (\ref{EQD1}) for $n=1$ is $$y(z)=\frac{1}{a-ib-z}$$ where $a-ib$ is an arbitrary complex constant. However, since $z\mapsto y(z)$ is analytic in $H^+$  we have  necessarily $b\geq 0.$ If $b>0$ one gets the Stieltjes transform of the Cauchy distribution $c_{a+ib}$, if $b=0$, then $\alpha=\delta_a.$ 

For $n=2$ things are more involved. Any solution of the differential equation $y''=2yy',$ corresponding to (\ref{EQD1}) for $n=2,$ which is analytic in $H^+$ satisfies $y'=y^2-C^2$ where $C$ is some complex constant. If $C=0$ we get that $y=\frac{1}{a-ib-z}$ as in the case $n=1$. In this case $\alpha$ is Cauchy or Dirac. Let us show now that taking $C\neq 0$ does not lead to an acceptable solution. We write first
$$1=\frac{y'}{y^2-C^2}=\frac{1}{2C}\left(\frac{y'}{y-C}-\frac{y'}{y+C}\right)$$ leading with an arbitrary constant $z_0$ to $y(z)=C\mathrm{cotanh}\; C(z_0-z).$ If $\Re C\neq 0$ the merophorphic function $z\mapsto \mathrm{cotanh} C(z_0-z)$ has poles in $H^+$ and $y$ would not be holomorphic in $H^+.$ If $C=ir$ is purely imaginary,  we observe that $y(z)=C\mathrm{cotanh}\; C(z_0-z)$
cannot be a Stieltjes transform since the condition $\lim_{t\rightarrow\pm \infty} y(z+t)=0$ is not fulfilled, the function $t\mapsto y(z+t)$ being periodic. 

Finally we consider the $d$-dimensional case. If $\alpha\in FT_d$ and if $\mu(n\alpha)=\alpha,$ let $f\in \R^d$ and denote by $\alpha_f$ the image of $\alpha$ by $x\mapsto \<f,x\>.$ Then $\mu(n\alpha_f)=\alpha_f$. If $n=1$ or $n=2$ we have seen that $\alpha_f$ is Cauchy: the definition of a Cauchy distribution in $\R^d$ implies the result. $\square$

\vspace{4mm}\noindent In the sequel, all the characterizations of the Cauchy distribution in $\R$ are extendable to $\R^d$ as done in Proposition 6.1, so we shall not mention it anymore and set $d=1$ from now on.

\vspace{4mm}\noindent\textbf{Proposition 6.2.} Let $\alpha\in FT_1$. Let $n<m$ any positive integers. Suppose that $\mu(n\alpha)=\mu(m\alpha)$ and let $y(z)=\inti \frac{\mu(n\alpha)(dw)}{w-z}.$ Then 
\begin{equation}\label{EQU4}\left(\frac{y^{(n-1)}}{(n-1)!}\right)^m=\left(\frac{y^{(m-1)}}{(m-1)!}\right)^n
\end{equation}
In particular if $m=n+1$ or if $m=n+2$ then $\alpha$ is Cauchy or Dirac.

\vspace{4mm}\noindent\textbf{Proof.} As usual we write $g(z)=-\inti\log (w-z)\alpha(dw).$ From Proposition 3.1 we have 
$$e^{ng(z)}=\inti \frac{\mu(n\alpha)(dw)}{(w-z)^n}=\frac{y^{(n-1)}(z)}{(n-1)!}$$From this (\ref{EQU4}) is plain.

Suppose now 
that $m=n+1$ and denote $Y=y^{(n-1)}/(n-1)!.$ From (\ref{EQU4}) we get 
$ (\frac{Y'}{n})^n=Y^{n+1}.$ Clearly $Y$ is not identically zero, since the  Stieltjes transform of a probability cannot be a polynomial. Select an open ball $U\subset H^+$ where $Y(z)\neq 0$ for all $z\in U.$ Therefore there exists a $n$th root of unity $\omega$ such that $Y'=n\omega Y^{1+\frac{1}{n}}.$ Integrating this differential equation we get that there exists  a complex number $a-ib$ such that $Y^{-1/n}=\omega(a-ib-z)$  leading to 
$\frac{y^{(n-1)}}{(n-1)!}=\frac{1}{(a-ib-z)^n}.$ Integrating $n-1$ times we get
$y(z)=P(z)+\frac{1}{a-ib-z}$ where $P$ is a polynomial with degree $<n.$ This is correct for $z\in U,$ but by analytic continuation it extends to the whole $H^+.$ Since $y$ is a Stieltjes transform $P=0$ and one concludes as the usual way that $b\geq 0$ and that  $\mu(n\alpha)$ is either Cauchy $c_{a+ib}$ or Dirac $\delta_a$ (from the Stieltjes transform of the Cauchy distribution). Since, again by Proposition 3.1, the map  $\alpha \mapsto \mu(n\alpha)$ is injective and from Corollary 2.3 $\mu(nc_{a+ib})=c_{a+ib}$ and $\mu(n\delta_{a})=\delta_a$ we conclude that $\mu(n\alpha)=\alpha$, so $\alpha$ is Cauchy or Dirac. 

Consider now the case $m=n+2$. From (\ref{EQU4}) we get $$\left(\frac{y^{(n-1)}}{(n-1)!}\right)^{n+2}=\left(\frac{y^{(n+1)}}{(n+1)!}\right)^n
$$ Again taking $Y=y^{(n-1)}/(n-1)!$ we get  $Y''=n(n+1)y^{(n+1)}/(n+1)!$ and finally 
$$\left(\frac{Y''}{n(n+1)}\right)^n=Y^{n+2}.$$ Using again a ball $U\subset H^+$ on which $Y(z)\neq 0$ there exists a $n$ th root of unity $\omega$ such that
$$Y''=n(n+1)\omega Y^{1+\frac{2}{n}}.$$ We now use a classical trick for ordinary differential equations of the form $Y''=f(Y',Y).$ From the implicit function theorem in the  analytic case, there exists an open set $V\subset U$ such that $z\mapsto Y(z)$ is injective  while restricted to $V$ and such that $Y(V)$ is open. As a consequence there exists an analytic function $p$ on $Y(V)$ such that $Y'(z)=p(Y(z))$ for $z\in V.$ Deriving we get $Y''(z)=p'(Y(z))p(Y(z))$ leading to 
$$2p'(Y(z))p(Y(z))=2n(n+1)\omega Y^{1+\frac{2}{n}}(z).$$ Thus integrating this differential equation in $p$ there exists a complex constant $C$ such that
$$p(Y(z))^2=(Y'(z))^2=n^2\omega(Y^{\frac{2n+2}{n}}(z)-C^{\frac{2n+2}{n}}).$$Now 
$Y(z)=\inti \frac{\alpha(dw)}{(w-z)^n}$ and $Y'(z)=n\inti \frac{\alpha(dw)}{(w-z)^{n+1}}$ imply that $C=0$ 
and that for some $2n$-th root of unity $\omega_1$  we have, for $z$ in some non empty open subset $V_1$ of $V$ 
$$Y'(z)=n\omega_1Y^{\frac{n+1}{n}}(z) $$ leading to the existence of a complex number $a-ib$ such that $Y^{-1/n}=\omega_1(z-a+ib).$ Since $\omega_1^{2n}=1$ we get $\omega_1^n=\pm 1$ and 
$$Y(z)=\pm \frac{1}{(a-ib-z)^n}.$$ Finally we get that $y(z)=P(z)\pm \frac{1}{a-ib-z}$ where $P$ is a polynomial. The fact that $y$ is a Stieltjes transform leads easily to $P=0$ and to $y(z)= \frac{1}{a-ib-z}$ where $b\geq 0$: this implies again that $\alpha$ is Cauchy or Dirac. $\square$

\vspace{4mm}\noindent\textbf{Proposition 6.3.} Let $\alpha\in FT_1.$ Let $N$ be an integer and suppose that $\mu(n\alpha)=\alpha$ for all $n\geq N.$ Then $\alpha$ is Cauchy or Dirac. 

\vspace{4mm}\noindent\textbf{Proof.} By Proposition 6.1, the hypothesis implies that for all $n\geq N$ we have 
$$y\frac{y^{(n-1)}}{(n-1)!}=\frac{y^{(n)}}{n!}$$ where $y(z)=\inti\frac{\alpha(dw)}{w-z}$ is the Stieltjes transform of $\alpha$, which is analytic in $H^+=\{z\in \mathbb{C}; \Im z>0\}.$ Since the above equality is true for all $n\geq N$ we deduce from it that for all $n\geq N$ we have 
\begin{equation}\label{EQU5}y^{n-N+1}\frac{y^{(N-1)}}{(N-1)!}=\frac{y^{(n)}}{n!}\end{equation}
Since $y$ is analytic in $H^+$ , when $z\in H^+$ the Taylor expansion of $t\mapsto y(z+t)$ converges for $|t|<\Im z$ and we can write for such $(z,t)$ 
\begin{eqnarray}\nonumber y(z+t)&=&\sum_{n=0}^{N-1}\frac{y^{(n)}(z)t^n}{n!}+\sum_{n=N}^{\infty}\frac{y^{(n)}(z)t^n}{n!}\\
\label{EQU52}&=&\sum_{n=0}^{N-1}\frac{y^{(n)}(z)t^n}{n!}+\frac{y^{(N-1)}(z)}{(N-1)!}\sum_{n=N}^{\infty}y^{n-N+1}(z)t^n\\
\label{EQU53}&=&\sum_{n=0}^{N-1}\frac{y^{(n)}(z)t^n}{n!}+\frac{y^{(N-1)}(z)}{(N-1)!}\frac{y(z)t^N}{1-ty(z)}
\end{eqnarray}
where (\ref{EQU52}) comes from (\ref{EQU5}). From (\ref{EQU53}) we get that $t\mapsto y(z+t)$  is a rational function. Since $y$ is analytic on $H^+$ this implies that (\ref{EQU53}) holds for all $z\in H^+$ and all real $t.$ 
We deduce from (\ref{EQU53}) by expanding the rational function $t\mapsto y(z+t)$ in partial fractions that there exists a polynomial $t\mapsto A_z(t)$ whose coefficients depend on $z$ such that 
\begin{equation}\label{EQU54}y(z+t)=A_z(t)+\frac{B_z}{1-ty(z)}\end{equation}where $B_z= \frac{y^{(N-1)}(z)}{(N-1)!}y(z)^{1-N}$ if $y(z)\neq 0$ and $B_z=0$ if $y(z)=0.$ The trick is now to observe that since $y$ is the  Stieltjes transform of the probability $\alpha$ we can write
$$\lim_{t\rightarrow \infty}ty(z+t)=\lim_{t\rightarrow \infty}t\inti\frac{\alpha(dw)}{w-z-t}=-1.$$ Applying this remark to  (\ref{EQU54}) we obtain that $A_z=0$, that $B_z=y(z)$ and finally that $y(z+t)=\frac{y(z)}{1-ty(z)}.$ Deriving with respect to  $t$ and setting $t=0$ we get $y'(z)=y^2(z),$ from which one concludes as in Proposition 6.1.$\square$

\vspace{4mm}\noindent\textbf{Proposition 6.4.} Let $\alpha\in FT_1$ and $0\leq b<c.$ Suppose that  $\nu=\mu(a\alpha)$  for all $a\in (b,c).$ Then  $\alpha=\nu$ is Cauchy or Dirac.

\vspace{4mm}\noindent\textbf{Proof.} Again with $g(z)=-\inti\log(w-z)\alpha(dw)$, with $z \in H^+,$ we can differentiate $n$ times with respect to $a \in (b,c)$ both sides of
$$\inti\frac{\nu(dw)}{(w-z)^a}=e^{ag(z)}.$$
We get for all $a\in (b,c)$
\begin{equation}\label{ET1}\inti [-\log(w-z)]^n\frac{\nu(dw)}{(w-z)^a}=e^{ag(z)}g(z)^n\end{equation}
The idea of the proof is to multiply both sides of (\ref{ET1}) by $t^n/n!$, to sum up in $n$, to invert sum and integral for finally getting 
$$\inti\frac{\alpha(dw)}{(w-z)^{a+t}}=e^{(a+t)g(z)}.$$
However the inversion of the sum and the integral needs some care. For this reason 
denote $u_n(w)=|-\log(w-z)|^n\frac{1}{|w-z|^a}$ and observe that $F(w,t)=\sum _{n=0}^{\infty}u_n(w)\frac{t^n}{n!}<\infty.$ If $0\leq t\leq a$ let us observe that  $$\inti F(w,t)\nu(dw)<\infty.$$ This obtained since $u_n(w)\leq (|\log|w-z||+\pi)^n\frac{1}{|w-z|^a}$ and therefore if $|w-z|>1$
$$F(w,t)\leq \frac{1}{|w-z|^a}e^{t|\log|w-z||+t\pi}=\frac{1}{|w-z|^{a-t}}e^{\pi t}$$
We now write from (\ref{ET1}) and the dominated convergence theorem 
\begin{eqnarray*}\inti \frac{\nu(dw)}{(w-z)^{a+t}}&=&e^{(a+t)g(z)}=\sum_{n=0}^{\infty}e^{ag(z)}g(z)^n\frac{t^n}{n!}\\&=&\sum_{n=0}^{\infty}\inti [-\log(w-z)]^n\frac{t^n}{n!}\frac{\alpha(dw)}{(w-z)^a}\\&=&\inti \frac{\alpha(dw)}{(w-z)^{a+t}}.\end{eqnarray*}
As a result $\alpha=\nu$ and furthermore $\mu((a+t)\alpha)=\alpha$ for all $t\in (0,a).$ By induction, we get easily that $\mu((a+t)\alpha)=\alpha$ for all $t>0.$ Now we apply Proposition 6.3 since $\mu(n\alpha)=\alpha$ for all integers $n$ large enough and the proof is complete. $\square$

\vspace{4mm}\noindent\textbf{Corollary 6.5.} If for a fixed $b$ and $c$ such that $0\leq b<c$ we have $\mu(b\alpha)=\mu(c\alpha)$ and if $\alpha$ has a mean then $\alpha$ is Dirac.

\vspace{4mm}\noindent\textbf{Proof.} If $b<a<c$ from Theorem 1.1 we have $\mu(c\alpha)\prec\mu(a\alpha)\prec \mu(b\alpha)$. From Comment 4 in Section 2 and from the hypothesis of the present corollary we have $\mu(c\alpha)=\mu(a\alpha)= \mu(b\alpha).$ Therefore  the hypothesis of  Proposition 6.4 is fulfilled and  $\alpha$ is Cauchy or Dirac. By since $\alpha$ has a mean, the first possibility is ruled out. $\square$

\vspace{4mm}\noindent\textbf{Proposition 6.6.} There exists a probability $\alpha \in FT_1$ such that  $\mu(\alpha) \notin FT_1.$

\vspace{4mm}\noindent\textbf{Proof.}
 Let us fix $1<a\leq 2$ and consider  
$$\alpha(dw)=\frac{a}{(1+\log(1+w))^{a+1}}\textbf{1}_{(0,\infty)}(w)\frac{dw}{1+w}.$$
With this definition, if $B\sim \alpha$, then $\Pr(\log(1+B)>t)=\frac{1}{(1+t)^a}$ for $t>0,$
so $\E(\log (1+B))<\infty.$ Let us compute 
\begin{eqnarray*}g(x)&=&-\int_{0}^{\infty}\log|x-w|\alpha(dw)=-\int_{0}^{\infty}\log|x-w|\frac{a}{(1+\log(1+w))^{a+1}}\frac{dw}{1+w}\\&=&-a\int_{0}^{\infty}\log|x+1-e^y|\frac{dy}{(1+y)^{a+1}}\\
g(e^u-1)&=&-a\int_{0}^{\infty}\log|e^u-e^y|\frac{dy}{(1+y)^{a+1}}=-u-a\int_{0}^{\infty}\log|1-e^{y-u}|\frac{dy}{(1+y)^{a+1}}
\end{eqnarray*}
From Cifarelli and Regazzini (1990) the density $f(x)$ of $X\sim \mu(\alpha)$ is, for $x>0,$ 
$$f(x)=\frac{1}{\pi}\sin\left(\pi\int_x^{\infty}\alpha(dw)\right)e^{g(x)}\sim_{x\rightarrow \infty}\left(\int_x^{\infty}\alpha(dw)\right)e^{g(x)}.$$ From this remark, $\E(\log(1+X))=\infty$  if and only if the integral $$I=\int_{0}^{\infty}\log (1+x)\left(\int_x^{\infty}\alpha(dw)\right)e^{g(x)}dx$$ diverges. Doing  in $I$ the change of variable $x=e^{u}-1$ we obtain
$$I=\int_{0}^{\infty}\frac{u}{(1+u)^a}e^{g(e^u-1)+u}du$$
From dominated convergence we have 
$$g(e^u-1)+u=-a\int_{0}^{\infty}\log|1-e^{y-u}|\frac{dy}{(1+y)^{a+1}}\rightarrow_{u\rightarrow \infty}0$$
Therefore $I$ diverges like the integral $J=\int_{0}^{\infty}\frac{udu}{(1+u)^{a}}$ since $1<a\leq 2.$ $\square$

\vspace{4mm}\noindent\textbf{Proposition 6.7.}  For $\alpha \in FT_1$ let $\mu_1(\alpha)=\mu(\alpha),$ and define by induction $\mu_n(\alpha)=\mu(\mu_{n-1}(\alpha)),$ if $\mu_{n-1}(\alpha)\in FT_1.$ Let $n\geq 2$ be an integer, and suppose that   $\alpha\in FT_1$ and $\mu_{k}(\alpha) \in FT_1$ for $k=2,\ldots,n-1$ and  $\mu_n(\alpha)=\alpha.$ Denote $y_j(z)=\inti \frac{\mu_j(\alpha)(dw)}{w-z},$ for $j=1,\ldots,n.$ Then
\begin{equation}\label{EQU7} (y'_1,\ldots,y'_n)=(y_ny_1,y_1y_2,\ldots,y_{n-1}y_n).\end{equation} In particular if $n=2$ then $\alpha$ is Cauchy or Dirac.

\vspace{4mm}\noindent\textbf{Proof.} With the convention $\mu_0(\alpha)=\alpha$ and the assumption $\mu_n(\alpha)=\alpha$ we can write for $j=1,\ldots,n:$
\begin{equation}\label{EQU72}\inti\frac{\mu_j(\alpha)(dw)}{w-z}=e^{g_{j-1}(z)}\end{equation} where $g_j(z)=-\inti\log (w-z)\mu_j(\alpha)(dw),$ for $j=0,\ldots,n-1.$ Since $g'_j=y_j,$ taking derivatives in (\ref{EQU72}) we get $y'_j=e^{g_{j-1}}g'_{j-1}=y_jy_{j-1},$ which is (\ref{EQU7}). If $n=2$ the differential system (\ref{EQU7}) gives $y'_1=y_1y_2=y_2'.$ Therefore there exists a complex constant $C$ such that $y_2=y_1+C$. If $C=0$ we  get  $y'_1=y_1^2$
leading to $\alpha$ being Cauchy and Dirac in the usual way. We are going to prove that $C\neq 0$ is impossible. Suppose the contrary: then, being  $y'_1=y_1(y_1+C)$ we get
$$\frac{1}{C}\left(\frac{y'_1}{y_1}-\frac{y'_1}{y_1+C}\right)=1$$ from which there exists a complex constant $z_0$ such that 
$y_1=\frac{C }{e^{-C(z-z_0)}-1}.$ The constant $z_0 $ cannot belong to $H^+:$  otherwise it is a pole of $y_1,$ which is impossible. Finally, if $\Re C\neq 0$ the function $y_1$ has poles in $H^+$, whereas if $C=ir$ is purely imaginary  the function $t\mapsto y_1(z+t)$  is periodic and this contradicts the fact that $y_1$ is  a Stieltjes transform. The proof is finished. $\square$

\section{Acknowledgments}  G.L thanks  Universit\`a di Roma  La Sapienza  and M.P. thanks  Universit\'e Paul Sabatier in Toulouse for their generous support during the preparation of this paper.

\section{References}\begin{enumerate}

\item \textsc{Cifarelli, D.M. and Melilli, E.} (2000) 'Some new results for Dirichlet priors' \textit{Ann. Statist.} \textbf{28}, 1390-1413.

\item \textsc{Cifarelli, D.M. and Regazzini, E.} (1990) 'Distribution functions of
means of a Dirichlet process'. \textit{Ann. Statist.}, 18, 429-442 (Correction in
\textit{Ann. Statist.} (1994) \textbf{22}, 1633-1634).

\item \textsc{Chamayou, J.-F.} (2000) Private communication. 

\item \textsc{Chamayou, J.-F. and  Letac, G.} (1991)
 'Explicit stationary distributions for compositions of random functions and product of random matrices.'
 {\it J. Theoret. Probab.} {\bf
4}, 3-36.
\item \textsc{Diaconis, P. and Kemperman, J.} (1996) 'Some new tools for Dirichlet priors' \textit{Bayesian Statistics 5 }, 97-106. Bernardo, Berger, Dawid and Smith (Eds), Oxford University Press.

\item \textsc{Diaconis, P. and Freedman, D.} (1999) 'Iterated random functions' \textit{SIAM Rev. }\textbf{41}, 45-76.

\item \textsc{Doss, H. and Sellke, T.} (1982)
 'The tails of probabilities chosen from a Dirichlet prior.'
 {\it Ann. Statist.} {\bf 10}, 1302-1305.

\item \textsc{Feigin, P. and Tweedie, R.L.} (1989)
 'Linear functionals and Markov chains associated with the Dirichlet processes.'
 {\it Math. Proc. Cambridge Philos. Soc.} {\bf 105}, 579-585.

\item \textsc{Ferguson, T.S.} (1973)
 'A Bayesian analysis of some nonparametric problems.'
 {\it Ann. Statist.} {\bf 1}, 209-230.

\item \textsc{Hannum, R.C., Hollander, M. and Landberg, N.A} (1981)
 'Random functionals of a Dirichlet process.' {\it Ann. Probab.} {\bf 9}, 665-670.
 
 \item \textsc{Hjort, N. L. and Ongaro, A.} (2005)
 'Exact inference for random Dirichlet means'
 {\it Stat. Inference  Stoch. Process.} {\bf 8}, 227-254.
 
 \item \textsc{James L.F.,  Lijoi A. and  Pr\"{u}nster} I. (2010) 'On the posterior distribution of classes of random means' \textit{Bernoulli} \textbf{16}, 155-180.

\item \textsc{James, L.} (2005)
 'Functionals of Dirichlet processes, the Cifarelli-Regazzini identity and Beta-Gamma processes.'
 {\it Ann. Statist.} {\bf 33}, 647-660.

\item \textsc{ Letac, G. and Piccioni, M.} (2014) 'Dirichlet random walks.'  To appear in \textit{J.  Appl. Probab.}, Dec. 2014.

\item \textsc{Lijoi, A. and Regazzini, E.} (2004) 'Means of a Dirichlet process and multiple hypergeometric functions' 
\textit{Ann. Probab.} \textbf{32}, 1469-1495.

\item \textsc{Lijoi, A. and Pr\"{u}nster, I.} (2009) 'Distributional properties of means of random probability measures' 
\textit{Statist. Surveys} \textbf{3}, 47-95.

\item \textsc{Muller, A. and Stoyan, D.} (2002)
\textit{Comparison methods for stochastic models and risks.} 
Wiley, Chichester.

\item \textsc{Ramachandran, B. and Rao, C.R.} (1970) 'Solutions of functional equations arising in some regression problems and a characterisation of the Cauchy laws' \textit{Sankhy\=a} Series A  \textbf{32}, 1-30.

\item \textsc{Samorodnitsky, G. and Taqqu, M.S.} (1994) {\it Stable non-Gaussian random processes.} Chapman and Hall, New York.

\item \textsc{Sethuraman, J.} (1994)
 'A constructive definition of Dirichlet priors.'
 {\it Statist. Sinica} {\bf 4}, 639-650.
 
 \item \textsc{Strassen, V.} (1965)
 'The existence of probability measures with given marginals.'
 {\it Ann. Math. Statist.} {\bf 36}, 423-439.

 \item \textsc{Wilks, S.} (1962) {\it Statistics.} Wiley, New York.

\item \textsc{Yamato, H.} (1984)
 'Characteristic functions of means of distributions chosen from a Dirichlet process.'
 {\it Ann. Probab.} {\bf 12}, 262-267.

\end{enumerate}

\end{document}